\documentclass[letterpaper,11pt,leqno]{article}
\usepackage{amsmath}
\usepackage{amsfonts}
\usepackage{amssymb}
\usepackage{amsthm}
\usepackage{array}
\usepackage{graphicx}
\usepackage{epsf}
\usepackage{url}

\newcommand{\x}{\mathbf{x}}
\newcommand{\dx}{\,d\x}
\newcommand{\deb}{\begin{eqnarray}}
\newcommand{\fin}{\end{eqnarray}}

\newtheorem{theorem}{Theorem}
\newtheorem{lemma}{Lemma}
\newtheorem{corollary}{Corollary}

\begin{document}

\setlength{\parindent}{0pt}

\begin{titlepage}

\begin{center}
{{\LARGE Some isoperimetric inequalities with}}
\vspace{5pt}

{{\LARGE application to the Stekloff problem}}

\vspace{26pt}

by
\end{center}

\vspace{1pt}

\begin{flushleft}
A. Henrot, Institut \'Elie Cartan, UMR7502 Nancy Universit\'e - CNRS - INRIA,
France,\\ e-mail : \url{antoine.henrot@iecn.u-nancy.fr}.\\
\vspace{.1in}

G.A. Philippin, D\'epartement de math\'ematiques et de statistique,
Universit\'e Laval, Sainte-Foy, Qu\'ebec, Canada, G1K 7P4, \\ e-mail : \url{gphilip@mat.ulaval.ca}.\\
\vspace{.1in}

A. Safoui, D\'epartement de math\'ematiques, Facult\'e des sciences --
Semlalia, Universit\'e Cadi Ayyad, Marrakech, Morocco, e-mail : \url{safoui@ucam.ac.ma}.\\

\vfill

\textbf{{Abstract}.} In this paper we establish isoperimetric
inequalities for the product of some moments of inertia. As an
application, we obtain an isoperimetric inequality for the product
of the $N$ first nonzero eigenvalues of the Stekloff problem in
$\mathbb{R}^N$.

\vfill


\vspace{22pt}


\end{flushleft}

\end{titlepage}

\section{{Introduction}}

Let $\x:=(x_1,...,x_N)$ be cartesian coordinates in $\mathbb{R}^N$,
$N\geq2$.  Let $J_k(\Omega)$ be the moment of inertia of $\Omega
\subset \mathbb{R}^N$ with respect to the plane $x_k=0$, defined as
\begin{eqnarray}
J_k(\Omega) := \int_{\Omega} {x_k}^2\dx\ ,\ \ k=1,...,N\ .
\end{eqnarray}

By summation over $k$, we obtain the polar moment of inertia of
$\Omega$ with respect to the origin denoted by $J_0(\Omega)$,
\begin{eqnarray}
J_0(\Omega) := \sum_{k=1}^{N} J_k(\Omega)\ .
\end{eqnarray}

Clearly, $J_0(\Omega)$ depends of the position of the origin. In
fact, $J_0(\Omega)$ is smallest when the origin coincides with the
center of mass of $\Omega$, i.e. when we have
\begin{eqnarray}
\int_{\Omega} x_k\dx=0\ ,\ \ k=1,...,N\ .
\end{eqnarray}

The following isoperimetric property is well known \cite{3,12} :

\vspace{.15in}
\begin{theorem}\label{theo1}
Among all domains $\Omega$ of prescribed $N$-volume, the ball
$\Omega^*$ centered at the origin has the smallest polar moment of
inertia, i.e. we have the isoperimetric inequality
\begin{eqnarray}\label{1.4}
J_0(\Omega) \geq J_0(\Omega^*)\ ,\ \ \Omega\in\mathcal{O}\ ,
\end{eqnarray}
where $\mathcal{O}$ is the class of all bounded domains of
prescribed $N$-volume.
\end{theorem}
\vspace{.15in}

Let $I_k(\Omega)$ be the moment of inertia of $\partial{\Omega}$ with respect to the plane $x_k=0$ defined as
\begin{eqnarray}
I_k(\Omega) := \int_{\partial\Omega} {x_k}^2\,ds\ ,\ \ k=1,2,...,N\ ,
\end{eqnarray}
where $ds$ is the area element of the boundary $\partial{\Omega}$ of $\Omega$.
By summation over $k$, we obtain the polar moment of inertia $I_0(\Omega)$ of $\partial{\Omega}$ with respect to the origin
\begin{eqnarray}
I_0(\Omega) := \sum_{k=1}^N I_k(\Omega)\ .
\end{eqnarray}
$I_0(\Omega)$ is smallest when the origin coincides with the center of mass of $\partial{\Omega}$, i.e. when we have
\begin{eqnarray}
\int_{\partial\Omega} {x_k}\,ds=0\ ,\ \ k=1,...,N\ .
\end{eqnarray}
Betta et al.\ \cite{1} have derived the following isoperimetric property.

\vspace{.1in}
\begin{theorem}
\begin{eqnarray}
I_0(\Omega)\geq I_0(\Omega^*)\ ,\ \ \Omega\in\mathcal{O}\ ,
\end{eqnarray}
with equality if and only if $\Omega$ coincides with $\Omega^*$.
\end{theorem}
\vspace{.1in}

The quantities of interest in this paper are
\begin{eqnarray}
J(\Omega):=\prod_{k=1}^N J_k(\Omega)\ ,
\end{eqnarray}
and
\begin{eqnarray}
I(\Omega):=\prod_{k=1}^N I_k(\Omega)\ .
\end{eqnarray}
We have several motivations to look at this products. Our main
motivation was to extend a classical result of Hersch, Payne and
Schiffer, see Theorem \ref{theoStek} below.

Another clear motivation is the following. Let us consider the
inertia matrix of a body $\Omega$ (we can do the same for its
boundary $\partial\Omega$), i.e. the matrix $\mathbf{M}$ whose
general term is
$$\mathbf{M}_{ij}=\int_{\Omega} {x_i}{x_j}\dx\ ,\ \ i,j=1,...,N\
.$$ The most classical invariants of this matrix are its trace and
its determinant and it is natural (for example for a mechanical
point of view) to ask for the domains which minimize these
invariants among all domains of prescribed $N$-volume. For the
trace, the answer is given in the Theorem \ref{theo1}. For the
determinant, we will see in section 2, that the ellipsoids symmetric
with respect to the planes $x_k=0$ are the minimizers of the
determinant: it will be a simple consequence of the study of the
product $J$.

A last motivation can be found in a paper of G. Polya, \cite{pol}.
Indeed, in this paper, the author was able to get the following
upper bound for the {\it torsional rigidity} $P(\Omega)$ of an
elastic beam with cross section $\Omega$:
$$J_0(\Omega)P(\Omega)\leq J(\Omega)$$
and, then was led to look for the minimizers of $J(\Omega)$ among
plane domains of given area.

\vspace{.1in} Let us now describe the content of this paper. In
Section 2 we establish the following isoperimetric property.
\vspace{.1in}
\begin{theorem}
\begin{eqnarray} \label{thm1}
J(\Omega) \geq J(E)\ ,
\end{eqnarray}
valid for all domains $\Omega\in\mathcal{O}$, with equality for all
ellipsoids $E\ (\in\mathcal{O})$ symmetric with respect to the
planes $x_k=0$, $k=1,...,N$.
\end{theorem}
We note that the isoperimetric inequality  (\ref{thm1}) also follows
from Blaschke's great contribution to affine geometry \cite{2}.
\vspace{.1in}

In Section 3 we establish the following isoperimetric property.

\vspace{.1in}
\begin{theorem}
\begin{eqnarray} \label{thm2}
I(\Omega) \geq I(\Omega^*)\ ,
\end{eqnarray}
valid for all convex domains $\Omega\in\mathcal{O}$, with equality
if and only if $\Omega=\Omega^*$.
\end{theorem}

\vspace{.1in}

As an application of (\ref{thm2}), we establish in Section 4 an
isoperimetric inequality for the product of the $N$ first nonzero
eigenvalues of the Stekloff problem in $\mathbb{R}^N$. This
inequality generalizes to dimension $N$ a previous two-dimensional
result of Hersch, Payne and Schiffer, see \cite{7}, \cite[Theorem
7.3.4]{8}:
\begin{theorem}\label{theoStek}
Let $\Omega$ be a convex domain in $\mathbb{R}^N$ and $0=p_1(\Omega)
< p_2(\Omega) \leq p_3(\Omega) \leq ...$ the eigenvalues of $\Omega$
for the Stekloff problem, see (\ref{4.1}). Then, the following
isoperimetric inequality holds
\begin{equation}\label{theost1}
\prod_{k=2}^{N+1}p_k(\Omega) \leq \prod_{k=2}^{N+1}p_k(\Omega^*)
\end{equation}
with equality if and only if $\Omega=\Omega^*$.
\end{theorem}
In (\ref{thm2}), the convexity of $\Omega$ may not be required. In
fact we show in Section 5 that (\ref{thm2}) remains valid for
nonconvex domains $\Omega$ in $\mathbb{R}^2$.


\section{{The functional $J$ and proof of (\ref{thm1})}}
\setcounter{equation}{0}

The proof of (\ref{thm1}) is based on the fact that $J(\Omega)$ is
not  affected by an affinity, i.e. when $\Omega$ is replaced by
\begin{eqnarray}
\Omega^{\prime}:=\{\ \x' := (t_1x_1,...,t_Nx_N) \in \mathbb{R}^N\ |\ \x\in\Omega\ \}\ ,
\end{eqnarray}
where $t_k$ are $N$ arbitrary positive constants such that
\begin{eqnarray}
\prod_{k=1}^N t_k = 1\ .
\end{eqnarray}

We compute indeed
\begin{eqnarray}\label{2.3}
J_k(\Omega)={t_k}^{-2}J_k(\Omega^{\prime})\ ,\ \ k=1,...,N\ ,
\end{eqnarray}
which implies
\begin{eqnarray}\label{2.4}
J(\Omega)=J(\Omega^{\prime})\ .
\end{eqnarray}

With the particular choice
\begin{eqnarray}
{t_k}^2:=(J(\Omega))^{\frac{1}{N}}J_k^{-1}(\Omega)\ ,\ \ k=1,...,N
\end{eqnarray}
in (\ref{2.3}), we obtain
\begin{eqnarray}\label{2.6}
J_k(\Omega^{\prime})=(J(\Omega))^{\frac{1}{N}}\ , \ \ k=1,...,N\ ,
\end{eqnarray}
i.e. the values of $J_k(\Omega^{\prime})$ are independent of $k$.
This shows that we have
\begin{eqnarray}\label{2.7}
\min_{\Omega\in\mathcal{O}} J(\Omega) =
\min_{\Omega^{\prime}\in\mathcal{O}^{\prime}} J(\Omega^{\prime})\ ,
\end{eqnarray}
where $\mathcal{O}^{\prime}$ ($\subset\mathcal{O}$) is the class of
all domains $\Omega^{\prime}$ of prescribed $N$-volume such that
$J_1(\Omega^{\prime})=...=J_N(\Omega^{\prime})$. Moreover we have by
(\ref{2.4}), (\ref{2.6}),
\begin{eqnarray}
J_0(\Omega^{\prime}):=\sum_{k=1}^N J_k(\Omega^{\prime})= N
(J(\Omega))^{\frac{1}{N}}=N (J(\Omega^{\prime}))^{\frac{1}{N}}\ ,
\end{eqnarray}
from which we obtain
\begin{eqnarray}\label{2.9}
J(\Omega^{\prime})=\left(\frac{1}{N}J_0(\Omega^{\prime})\right)^N\ .
\end{eqnarray}

Combining (\ref{2.7}) and (\ref{2.9}), and making use of
(\ref{1.4}), we are led to
\begin{equation}
\begin{array}{lrcl}
 (2.10) \ \ \ & \displaystyle \min_{\Omega \, \in \,
\mathcal{O}} J(\Omega) &=& \displaystyle \min_{\Omega^{\prime} \,
\in \,\mathcal{O}^{\prime}}
\left\{\frac{1}{N}J_0(\Omega^{\prime})\right\}^N
= \frac{1}{N^N}\left( \min_{\Omega^{\prime} \, \in \, \mathcal{O}^{\prime}} \left\{J_0(\Omega^{\prime})\right\} \right)^N\\
& &=& \displaystyle \frac{1}{N^N}\left(J_0(\Omega^*)\right)^N=J(\Omega^*)=J(E)\ , \nonumber
\end{array}
\end{equation}
which is the desired result.

\vspace{.1in} We give now an application to the minimization of the
determinant of the inertia matrix. We assume that the origin $O$ is
at the center of mass of the domains we consider.
\begin{corollary}
Let $\mathbf{M}(\Omega)$ be the inertia matrix of the domain
$\Omega$ i.e. the matrix whose general term is
$$\mathbf{M}_{ij}=\int_{\Omega} {x_i}{x_j}\dx\ ,\ \ i,j=1,...,N$$
and let $D(\Omega)$ be its determinant. Then,
\begin{equation}\label{2.11}
    D(\Omega)\geq D(E)\,
\end{equation}
valid for all domains $\Omega\in\mathcal{O}$, with equality for all
ellipsoids $E\ (\in\mathcal{O})$ symmetric with respect to the
planes $x_k=0$, $k=1,...,N$.
\end{corollary}
Indeed, since $\mathbf{M}(\Omega)$ is symmetric, there exists an
orthogonal matrix $T\in O^+(N)$ and a diagonal matrix $\Delta$ such
that $\mathbf{M}(\Omega)=T^T{\Delta}T$. Actually, $\Delta$ is the
inertia matrix of the domain $(T(\Omega)$ obtained from $\Omega$ by
some rotation. Now, the determinant being invariant through such a
similarity transformation, we have according to (\ref{thm1}):
\begin{equation}\label{2.12}
    D(\Omega)=D(T(\Omega))=J(T(\Omega))\geq J(E)=D(E)
\end{equation}
the last equality in (\ref{2.12}) coming from the fact that $E$ is
symmetric with respect to each plane of coordinates, so its inertia
matrix is diagonal. This proves the desired result.


\section{{The functional $I$ and proof of (\ref{thm2})}
\setcounter{equation}{0}}

We assume in this section that $\Omega$ is convex and choose the
origin at the center of mass of $\partial\Omega$. We introduce the
family of parallel domains
\begin{eqnarray}
\Omega_h&:=&\Omega+B_h(0)\\
&=&\{\ \mathbf{z}=\x+\mathbf{y} \in \mathbb{R}^N\ |\ \x\in\Omega\ ,\ \mathbf{y}\in B_h(0)\ \}\ ,\nonumber
\end{eqnarray}
where $B_h(0)$ is the $N$-ball of radius $h>0$ centered at  the
origin. It follows from the Brunn-Minkowski theory that the function
$f(h):=\left|\Omega_h\right|^{1/N}$ is concave with respect to the
parameter $h$. H. Minkowski made use of this basic property to
derive the famous classical isoperimetric geometric inequality for
convex bodies. Since our approach will be patterned after his
argument, we indicate here briefly Minkowski's method. The concavity
of $f(h)$ implies that $f^{\prime}(h)$ is a monotone decreasing
function. This leads to the inequality
\begin{eqnarray}
f^{\prime}(h)=\frac{1}{N}\left|\Omega_h\right|^{\frac{1-N}{N}}
\left|\Omega_h\right|^{\prime}=\frac{1}{N}\left|\Omega_h\right|^{\frac{1-N}{N}}\left|\partial\Omega_h\right|\geq C\ ,
\end{eqnarray}
where $C$ has to coincide with the value of $f^{\prime}(h)$ for a
ball since $\Omega_h$ approaches a large ball as $h$ increases to
infinity. Evaluated at $h=0$, (3.2) leads to the well known
isoperimetric geometric inequality
\begin{eqnarray}
\left|\partial\Omega\right|^N\geq N^N\omega_N\left|\Omega\right|^{N-1}\ ,
\end{eqnarray}
where $\omega_N:=\pi^{N/2}/\,\Gamma\left(\frac{N}{2}+1\right)$ is
the volume of the unit ball in $\mathbb{R}^N$. We refer the reader
to the basic books of Bonnesen-Fenchel \cite{3} and Hadwiger
\cite[6]{5} for details. The proof of (\ref{thm2}) makes use of the
following two lemmas.

\vspace{.15in}
\begin{lemma}
With the notations of Section 1, we have
\begin{eqnarray}
J_k(\Omega_h)=J_k(\Omega)+hI_k(\Omega)+\sum_{j=2}^{N+2}c_jh^j\ ,\ \ h\geq 0\ ,
\end{eqnarray}
where the coefficients $c_j$ are some geometric quantities associated to $\partial\Omega$.
\end{lemma}
\vspace{.15in}

As a consequence of (3.4), we have
\begin{eqnarray}
\left.{J_k}'(h)\right|_{h=0}=I_k(\Omega)\ .
\end{eqnarray}

For the proof of Lemma 1, we evaluate $\int_{\Omega_h\setminus\Omega} x_k^2\dx=J_k(\Omega_h)-J_k(\Omega).$
The computation of this integral will be easy if we introduce normal coordinates $\mathbf{s}:=(s_1,...,s_N)$ such that

\vspace{.15in}

$
\begin{array}{ll}
(3.6)\ \ &
\Omega_h\setminus\Omega=
\left\{
\x(\mathbf{s})=\mathbf{r}(s_1 , ... , s_{N-1})+s_N\mathbf{n}(s_1, ... ,s_{N-1})\ ,\ \
\mathbf{s}\in B\times [0,h]
\right\}\ ,
\end{array}
$

\vspace{.15in}

\setcounter{equation}{6}

where $\mathbf{r}(s_1, ... ,s_{N-1}),\ (s_1, ... ,s_{N-1})\in B$, is
a parametric representation of $\partial\Omega$, and
$\mathbf{n}(s_1,...,s_{N-1})$ is the unit normal vector of
$\partial\Omega$. In terms of the new variables $\mathbf{s}$, the
volume element of $\mathbb{R}^N$ is $dV=\Delta \, ds_1...ds_N$, with
\begin{eqnarray}
\Delta&=&\det \left|\frac{\partial x_k}{\partial s_j}\right|\\
&=&\det \left| \frac{\partial \mathbf{r}}{\partial s_1}+s_N
\frac{\partial \mathbf{n}}{\partial s_1} \, , \, ...\, , \,
\frac{\partial \mathbf{r}}{\partial s_{N-1}}+s_N \frac{\partial
\mathbf{n}}{\partial s_{N-1}} \, , \, \mathbf{n} \, \right|\ .\nonumber
\end{eqnarray}

$\Delta$ is obviously a polynomial in $s_N$ of degree $(N-1)$.
We have
\begin{eqnarray}
\Delta=\Delta_0+\sum_{j=1}^{N-1}\widetilde{c_j}\left(s_N\right)^j\ ,
\end{eqnarray}
where $\Delta_0$ is $\Delta$ evaluated on $\partial\Omega$. We then obtain
\begin{eqnarray}
\int_{\Omega_h\setminus\Omega} x_k^2\dx &=& \int_{B\times[0,h]}
(x_k^2 \Delta_0\,+\sum_{j=1}^{N-1}
x_k^2 \widetilde{c_j}\left(s_N )^j \right)\,ds_1\ldots ds_N \nonumber\\
&=&hI_k(\Omega)+\sum_{j=2}^{N+2} c_jh^j\ ,\ \ h\geq0\ ,\nonumber
\end{eqnarray}
which is the desired result (3.4).

\vspace{.15in}

The next lemma follows from an extension of the Brunn-Minkowski inequality established by H. Knothe \cite{10}.

\vspace{.15in}
\begin{lemma}
The functions $g(h):=\left( J_k\left(\Omega_h\right) \right)^{\frac{1}{N+2}}$ are concave for $h\geq 0$, $k=1,...,N$.
\end{lemma}
\vspace{.15in}

As a direct consequence of Lemmas 1 and 2, we have in analogy to (3.2)
\begin{eqnarray}\label{3.9}
g'(h)=\frac{1}{N+2}\left( J_k\left(\Omega_h\right) \right)^{-\frac{N+1}{N+2}}I_k(\Omega_h)\geq C\ ,\ \ k=1,...,N\ ,
\end{eqnarray}
where $C$ has to coincide with the value of $g'(h)$ for a $N$-ball.
Using $J_k(B_R)=\frac{R^{N+2}\omega_N}{N+2}$ for the ball of radius
$R$, we get $C=\left(\frac{\omega_N}{N+2}\right)^{1/(N+2)}$.
Evaluated at $h=0$, (\ref{3.9}) leads to the following isoperimetric
inequalities
\begin{eqnarray}\label{3.10}
(I_k(\Omega))^{N+2} \geq (N+2)^{N+1}\omega_N(J_k(\Omega))^{N+1}\ ,\
\ k=1,...,N\ ,
\end{eqnarray}
with equality if and only if $\Omega=\Omega^*$. Making use of
(\ref{3.10}) and (1.11), we obtain
\begin{equation}
\begin{array}{lrcl}
\ \ \ & (I(\Omega))^{N+2} &\geq& (N+2)^{N(N+1)}(\omega_N)^N(J(\Omega))^{N+1}\\
& &\geq&
(N+2)^{N(N+1)}(\omega_N)^N(J(\Omega^*))^{N+1}=(I(\Omega^*))^{N+2}\ ,
\nonumber
\end{array}
\end{equation}
\setcounter{equation}{12}
with equality if and only if $\Omega=\Omega^*$. This achieves the proof of inequality (1.12) for convex domains $\Omega$.


\section{{Application to the Stekloff problem}}
\setcounter{equation}{0}

In this section we consider the Stekloff eigenvalue problem defined in a bounded convex domain $\Omega$ in $\mathbb{R}^N$, $N\geq 2$.
\begin{eqnarray}\label{4.1}
\left\{ \begin{array}{lcl} \displaystyle{\Delta u = 0}  &,& \displaystyle{\x:=(x_1,...,x_N) \in \mathbb{R}^N} \ ,\\
\displaystyle{\frac{\partial u}{\partial n} = p u } &,& \displaystyle{ \x \in \partial \Omega} \ .\end{array}
\right.
\end{eqnarray}
In (4.1), $\frac{\partial u}{\partial n}$ is the exterior normal
derivative of  $u$ on $\partial \Omega$. It is well known \cite{13}
that there are infinitely many eigenvalues $0=p_1 < p_2 \leq p_3
\leq ...$ for which (4.1) has nontrivial solutions, also called
eigenfunctions, and denoted by $u_1 (=\mbox{const}.),u_2,u_3,...$
Let $\Sigma_k$ be the class of functions defined as
\begin{eqnarray}
\Sigma_k := \left\{\ v\in H^1(\Omega)\ ,\ \int_{\partial\Omega} v u_j \, ds=0\ ,\ j=1,...,k-1\ \right\}\ ,
\end{eqnarray}
where $H^1(\Omega)$ is the Sobolev space of functions in
$L^2(\Omega)$  whose partial derivatives are in $L^2(\Omega)$. Let
$R[v]$ be the the Rayleigh quotient associated to the problem (4.1)
defined as
\begin{eqnarray}
R[v]:=\frac{\int_\Omega|\nabla v|^2\dx}{\int_{\partial\Omega} v^2 \, ds}\ .
\end{eqnarray}

It is well known that the eigenvalue $p_k$ has the following variational characterization \cite[8]{7}
\begin{eqnarray}
p_k=\min_{v\,\in\,\Sigma_k} R[v]\ ,\ \ k=2,3,4,...\ .
\end{eqnarray}

Unfortunately, (4.4) is of little practical use for estimating
$p_k$,  since it requires the  knowledge of the eigenfunctions
$u_j$, $j=1,...,k-1$. The following variational characterization due
to H. Poincar\'e overcomes this difficulty. Let $v_k (\not \equiv 0)
\in H^1(\Omega),\ k=1,...,n$ be $n$ linearly independent functions
satisfying the conditions
\begin{eqnarray}
\int_{\partial\Omega} v_k \, ds = 0\ ,\ \ k=1,...,n\ .
\end{eqnarray}
Let $L_n$ be the linear space generated by $v_k,\ k=1,...,n$. The
Rayleigh quotient of $v:=\sum_{k=1}^n c_k v_k$ is the ratio of two
quadratic forms of the $n$ variables $c_1,...,c_n$. We have
\begin{eqnarray}
R[v]:=\frac{\int_\Omega|\nabla v|^2\dx}{\int_{\partial\Omega} v^2 \, ds}=
\frac{\sum_{i,j=1}^n a_{ij} c_i c_j}{\sum_{i,j=1}^n b_{ij} c_i c_j}\ ,
\end{eqnarray}
with
\begin{eqnarray}
a_{ij}&:=&\int_{\Omega} \nabla v_i \nabla v_j\dx\ ,\\
b_{ij}&:=&\int_{\partial\Omega} v_i v_j \, ds \ .
\end{eqnarray}

Note that the matrices $A:=(a_{ij}), B:=(b_{ij})$ are positive
definite.  Let $0  < p_2' \leq p_3' \leq ...\leq p_{n+1}' $ be the
$n$ roots of the characteristic equation
\begin{eqnarray}
\det|A-p B|=0\ .
\end{eqnarray}

Poincar\'e's variational principle \cite{11} asserts that
\begin{eqnarray}
p_k \leq p_k' \ ,\ \ k=2,...,n+1\ .
\end{eqnarray}

By means of a translation followed by an appropriate rotation, the following conditions will we satisfied
\begin{eqnarray}
\int_{\partial \Omega} x_k\,ds&=&0 \ , \ \ k=1,...,N\ ,\\
\int_{\partial \Omega} x_k x_j \,ds&=&0 \ , \ \ k \neq j\ .
\end{eqnarray}

The $N$ functions defined as
\deb
v_k:=x_k(I_k(\Omega))^{-1/2}\ ,\ \ k=1,...,N
\fin
are admissible for the Poincar\'e principle. We then compute with the notation of Section 1
\deb
A&=&|\Omega|\,\mbox{diag}\big(\, I_1^{-1}(\Omega),...,I_N^{-1}(\Omega) \, \big)\ ,\\
\nonumber {}\\
B&=&\mbox{diag}\big(\, 1,...,1 \, \big)\ .
\fin

The $N$ roots of the characteristic equation (4.9) are then
$|\Omega|I_k^{-1}(\Omega)$, \\$k=1,...,N$. We then obtain from (4.10) with $n=N$
\deb
\prod_{k=2}^{N+1}p_k(\Omega) &\leq& \prod_{k=2}^{N+1}p_k'=|\Omega|^N I^{-1}(\Omega) \\
 &\leq& |\Omega|^N I^{-1}(\Omega^{*}) = \frac{\omega_N}{|\Omega|}=\prod_{k=2}^{N+1}p_k(\Omega^*)\ .  \nonumber
\fin

In (4.16), we have used the isoperimetric inequality (1.12) and the
fact that $x_k$, $k=1,...,N$, are the $N$ first nonzero eigenvalues
of $\Omega^*$.

Note that (4.16) is an improvement (for convex $\Omega$!) of the
following inequality \deb \sum_{k=2}^{N+1}\frac{1}{p_k(\Omega)} \geq
\sum_{k=2}^{N+1} \frac{1}{p_k(\Omega^*)}\ , \fin obtained by Brock
\cite{4}.


\section{{Appendix}
\setcounter{equation}{0}}

Since our proof of (1.12) is valid only for convex $\Omega$, we
indicate  in this section a proof inspired by an old paper of A.
Hurwitz \cite{9} that does not require convexity of $\Omega \subset
\mathbb{R}^2$. Let $L$ be the length of $\partial\Omega$ and $s$ be
the arc length on $\partial\Omega$. Consider the following
parametric representation of $\partial\Omega$ : \deb
(x(\sigma),y(\sigma))\ ,\ \ \sigma:=\frac{2\pi}{L}s \ \ \in \
[0,2\pi]\ . \fin

Clearly $x(\sigma)$ and $y(\sigma)$ are $2\pi$-periodic functions of
$\sigma$ whose associated Fourier series are of the form
\deb
\left\{ \begin{array}{l} \displaystyle{ x(\sigma) = \frac{1}{2}a_0 +
\sum_{k=1}^{\infty}(a_k\cos k\sigma +{a'_k}\sin k\sigma ) }\ ,  \\
\displaystyle{ y(\sigma) = \frac{1}{2}b_0 +
\sum_{k=1}^{\infty}(b_k\cos k\sigma +{b'_k}\sin k\sigma ) } \ .
\end{array} \right.
\fin

The Fourier coefficients $a_k,\ {a'_k},\ b_k,\ {b'_k}$ have to be
determined in order to minimize $I(\Omega)=I_1(\Omega)I_2(\Omega)$.
>From (5.2) and Parseval's identity, we obtain \deb
\frac{L^2}{2\pi}&=&\frac{1}{\pi}\int_0^{2\pi} \left\{ \left( \frac{d
x}{d\sigma} \right)^2 + \left( \frac{d y}{d\sigma} \right)^2
\right\}\,d\sigma \\ &=& \sum_{k=1}^\infty k^2(a_k^2+{{a'_k}}^2 +b_k^2+{{b'_k}}^2)\ ,\nonumber \\
\nonumber{}\\
|\Omega|&=&\int_0^\pi x(\sigma) \frac{dy}{d\sigma}\, d\sigma = \pi \sum_{k=1}^\infty k ( a_k{b'_k}-{a'_k}b_k )\ ,\\
\nonumber{}\\
  I_1(\Omega)&=& \int_{\partial\Omega}x^2(s)\,ds=\frac{L}{2\pi}
  \int_0^{2\pi} x^2(\sigma)\, d\sigma=\frac{L}{2}\left\{ \frac{1}{2} {a_0}^2+a^2 \right\}\ ,\\
\nonumber{}\\
  I_2(\Omega)&=& \int_{\partial\Omega}y^2(s)\,ds=\frac{L}{2\pi}
  \int_0^{2\pi} y^2(\sigma)\, d\sigma=\frac{L}{2}\left\{ \frac{1}{2} {b_0}^2+b^2 \right\}\ ,
\fin with \deb a^2:= \sum_{k=1}^\infty ({a_k}^2+{{a'_k}}^2) ,\ \
b^2:=\sum_{k=1}^\infty ({b_k}^2+{{b'_k}}^2)\ . \fin Clearly, we must
choose $a_0=b_0=0$ since $L$ and $|\Omega|$ are  independent of
$a_0,\ b_0$. The other Fourier coefficients may be determined using
Lagrange's method, consisting in finding the critical points of the
Lagrange function defined as \deb
F(a,a',b,b')&:=& 4 I(\Omega) + \frac{\lambda}{\pi}|\Omega|\\
\nonumber &=& L^2 a^2 b^2 - \lambda \sum_{k=1}^\infty (
a_k{b'_k}-{a'_k}b_k )\ , \fin where $\lambda$ is a multiplier. This
leads to the following system of equations \deb
\frac{\partial F}{\partial a_k} &=& 4 \pi^2 k^2 a^2 b^2 a_k + 2L^2b^2a_k-\lambda k b'_k=0\ ,\\
\frac{\partial F}{\partial a'_k} &=& 4 \pi^2 k^2 a^2 b^2 {a'_k} + 2L^2b^2{a'_k}+\lambda k b_k=0\ ,\\
\frac{\partial F}{\partial b_k} &=& 4 \pi^2 k^2 a^2 b^2 b_k + 2L^2a^2b_k+\lambda k {a'_k}=0\ ,\\
\frac{\partial F}{\partial {b'_k}} &=& 4 \pi^2 k^2 a^2 b^2 {b'_k} +
2L^2a^2{b'_k}-\lambda k a_k=0\ . \fin

>From (5.9), (5.12), we obtain \deb a_k M_k = 0 \ , \fin with \deb
M_k := 4a^2b^2(2\pi^2k^2a^2+L^2)(2\pi^2k^2b^2+L^2)-k^2\lambda^2\ .
\fin

>From (5.10), (5.11), we obtain \deb {a'_k} M_k=0\ . \fin

We conclude that either $a_k={a'_k}=0$, $k=1,2,3,...$, which is
absurd, or that
\begin{equation}
\begin{array}{r}\lambda^2 = 4 a^2 b^2 k^{-2} (2\pi^2 k^2 a^2 + L^2)(2 \pi^2 k^2 b^2 + L^2) = 4a^2 b^2 f(k^2)\ ,\end{array}
\end{equation}
with
\deb
f(t):=\frac{1}{t}(2\pi^2 t a^2 + L^2)(2\pi^2 t b^2 + L^2)\ ,\ \ t \geq 1 \ .
\fin

Since $f(t)$ is convex for $t \geq 1$, equation (5.16) can be
satisfied  for at most two positive integers $k_1 \leq k_2$. We then
conclude that the parametric representation of $\partial \Omega^*$
is of the form

\vspace{.15in}

$
\begin{array}{ll}
(5.18) &
\left\{ \begin{array}{l} \displaystyle{ x(\sigma) = a_{k_1}\cos( k_1\sigma)
+a_{k_1}'\sin (k_1\sigma ) + a_{k_2}\cos( k_2\sigma) +a_{k_2}'\sin (k_2\sigma ) }\ ,  \\
{}\\
\displaystyle{ y(\sigma) = b_{k_1}\cos( k_1\sigma) +b_{k_1}'\sin (k_1\sigma )
+ b_{k_2}\cos( k_2\sigma) +b_{k_2}'\sin (k_2\sigma ) } \ .
\end{array} \right.
\end{array}
$

\vspace{.15in}

\setcounter{equation}{18}

Further restrictions on the Fourier coefficients $a_{k_j},\
a_{k_j}', \ b_{k_j},\ b_{k_j}'$, $j=1,2$, are imposed by the
condition \deb \left( \frac{dx}{d\sigma} \right)^2 + \left(
\frac{dy}{d\sigma} \right)^2=\left( \frac{L}{2\pi}
\right)^2=\mbox{const.} \fin

>From (5.18), we compute
\begin{equation}
\begin{array}{l}
\displaystyle{\left( \frac{dx}{d\sigma} \right)^2 + \left( \frac{dy}{d\sigma}
\right)^2 = c_0+c_1\cos(2k_1\sigma)+c_2\sin(2k_1\sigma)} \\
{}\\
 \displaystyle{+ c_3\cos(2k_2\sigma) + c_4\sin(2k_2\sigma) + c_5\cos(k_1-k_2)\sigma} \\
{}\\
 \displaystyle{+ c_6 \cos(k_1+k_2)\sigma + c_7\sin(k_1-k_2)\sigma + c_8\sin(k_1+k_2)\sigma}\ , \end{array}
\end{equation}


with
\begin{equation*}
\begin{array}{lrcl}
(5.21)\ \  & c_0 & := &\frac{1}{2} k_1^2 (a_{k_1}^2 + a_{k_1}^{'2} +
b_{k_1}^2 + b_{k_1}^{'2} ) + \frac{1}{2} k_2^2 (a_{k_2}^2 +
a_{k_2}^{'2} + b_{k_2}^2 + b_{k_2}^{'2} )\ ,\\
{}\\
(5.22) & c_1 & := &\frac{1}{2} k_1^2 (a_{k_1}^{'2} + b_{k_1}^{'2} - a_{k_1}^2 - b_{k_1}^2 )\ ,\\
{}\\
(5.23) & c_2 & := &- k_1^2 ({a_{k_1}}{a_{k_1}'} + {b_{k_1}}{b_{k_1}'} )\ ,\\
{}\\
(5.24) & c_3 & := &\frac{1}{2} k_2^2 (b_{k_2}^{'2} + a_{k_2}^{'2} - b_{k_2}^2 - a_{k_2}^2 )\ ,\\
{}\\
(5.25) & c_4 & := &- {k_2}^2 ({a_{k_2}}{a_{k_2}'} + {b_{k_2}}{b_{k_2}'} )\ ,\\
{}\\
(5.26) & c_5 & := & {k_1}{k_2} ({a_{k_1}}{a_{k_2}} + {b_{k_1}}{b_{k_2}} + {a_{k_1}'}{a_{k_2}'} + {b_{k_1}'}{b_{k_2}'} )\ ,\\
{}\\
(5.27) & c_6 & := & {k_1}{k_2} (-{a_{k_1}}{a_{k_2}} - {b_{k_1}}{b_{k_2}} + {a_{k_1}'}{a_{k_2}'} + {b_{k_1}'}{b_{k_2}'} )\ ,\\
{}\\
(5.28) & c_7 & := & {k_1}{k_2} ({a_{k_1}'}{a_{k_2}} + {b_{k_1}'}{b_{k_2}} - {a_{k_1}}{a_{k_2}'} - {b_{k_1}}{b_{k_2}'} )\ ,\\
{}\\
(5.29) & c_8 & := & -{k_1}{k_2} ({a_{k_1}}{a_{k_2}'} + {b_{k_1}}{b_{k_2}'} + {a_{k_1}'}{a_{k_2}} + {b_{k_1}'}{b_{k_2}} )\ .
\end{array}
\end{equation*}
\setcounter{equation}{29}

Suppose now that (5.19) is satisfied for two positive integers $k_1\neq k_2$. Then we must have
\deb
c_j=0\ ,\ \ j=1,...,8 \ \ \mbox{if}\ \ k_1 \neq 3k_2\ ,
\fin
or
\deb
c_1=c_2=c_6=c_8=c_3+c_5=c_4+c_7=0\ \ \mbox{if}\ \ k_1=3k_2\ .
\fin

A careful investigation of these two cases shows that we must have
either $a_{k_1}=a_{k_1}'=b_{k_1}=b_{k_1}'=0$, or
$a_{k_2}=a_{k_2}'=b_{k_2}=b_{k_2}'=0$. For the sake of brevity we
omit the computational details to confirm this assertion. In any
case the parametric representation of $\partial\Omega^*$ takes the
following form

\vspace{.15in}

$
\begin{array}{ll}
(5.32) &
\left\{ \begin{array}{l} \displaystyle{ x(\sigma) = a_{k_0}\cos k_0\sigma +a_{k_0}'\sin k_0\sigma }\ ,  \\
{} \\
\displaystyle{ y(\sigma) = b_{k_0}\cos k_0\sigma +b_{k_0}'\sin k_0\sigma  \ ,\ \ \sigma \in [0,2\pi]}\ ,
\end{array} \right.
\end{array}
$

\vspace{.15in}

for some positive integer $k_0$. But since $\partial\Omega^*$ makes
only  one loop around the origin, we must actually have $k_0=1$. We
then obtain

\vspace{.15in}

$
\begin{array}{ll}
(5.33) &
\left\{ \begin{array}{l} \displaystyle{ x(\sigma) = a_{1}\cos\sigma +a_{1}'\sin \sigma }  \\
{} \\
\displaystyle{ y(\sigma) = b_{1}\cos\sigma +b_{1}'\sin \sigma}  \ ,
\end{array} \right.
\end{array}
$

\vspace{.15in}

with

\vspace{.15in}

$
\begin{array}{ll}
(5.34) &
a_1^2+b_1^2=a_1'^2+b_1'^2\ , \ \ a_1a_1'+b_1b_1'=0\ .
\end{array}
$

\vspace{.15in}

Finally (5.33), (5.34) lead to

\vspace{.15in}

$
\begin{array}{ll}
(5.35) &
x^2(\sigma)+y^2(\sigma)=a_1^2+b_1^2=\mbox{const.}
\end{array}
$

\vspace{.15in}

(5.35) shows that if there exists a minimizer $\widetilde{\Omega}$
of $I(\Omega)$, $\widetilde{\Omega}\in \mathcal{O}$, then it must be
a disc centered at the origin. The existence of $\widetilde{\Omega}$
(among other similar results) will be established in a forthcoming
paper of A. Henrot, \cite{hen}.

\vspace{.3in}

\textbf{Acknowledgment :} This research has been initiated at the AIM
workshop on low eigenvalues of Laplace and Schr\"odinger operators, Palo Alto (Ca.), may 2006.

\end{document}